\documentclass[12pt]{article}
\usepackage{amsmath}
\usepackage{amssymb}
\usepackage{amscd}
\usepackage{amsthm}
\usepackage{graphicx}
\usepackage{epsfig}

\numberwithin{equation}{section}

\theoremstyle{plain}
\newtheorem{theorem}[equation]{Theorem}
\newtheorem{thm}[equation]{Theorem}
\newtheorem{prop}[equation]{Proposition}

\newtheorem{lem}[equation]{Lemma}
\newtheorem{cor}[equation]{Corollary}

\newtheorem{rem}[equation]{Remark}

\theoremstyle{definition}

\newtheorem{defn}[equation]{Definition}

\newcommand{\R}{\mathbb R}

\newcommand{\C}{\mathbb C}
\newcommand{\Z}{\mathbb Z}
\newcommand{\Q}{\mathbb Q}

\renewcommand{\P}{\mathbb P}
\renewcommand{\H}{\mathbb H}

\newcommand{\Area}{\operatorname{Area}}
\newcommand{\Diff}{\operatorname{Diff}}
\newcommand{\Zero}{\operatorname{Zero}}
\newcommand{\Span}{\operatorname{Span}}
\newcommand{\Proj}{\operatorname{Proj}}
\newcommand{\Vol}{\operatorname{Vol}}
\newcommand{\Image}{\operatorname{Image}}

\newcommand{\acts}{\curvearrowright}

\newcommand{\al}{\alpha}

\def\Del{\Delta}

\def\ga{\gamma}
\def\Ga{\Gamma}

\def\La{\Lambda}
\def\si{\sigma}
\def\Si{\Sigma}

\newcommand{\ov}{\overrightarrow}

\def\om{\omega}
\def\Om{\Omega}

\def\si{\sigma}
\def\Si{\Sigma}

\def\be{\beta}
\def\del{\delta}

\def\t{\tilde}

\newcommand{\ol}{\overline}

 \newenvironment{dedication}
        {\vspace{6ex}\begin{quotation}\begin{center}\begin{em}}
        {\par\end{em}\end{center}\end{quotation}}

\begin{document}

\title{Periods of abelian differentials and dynamics}
\author{Misha Kapovich}
\date{\today}
\maketitle

\begin{abstract}
Given a closed oriented surface $S$ of genus $\ge 3$ 
we describe those cohomology classes $\chi\in H^1(S,\C)$ which appear 
as the period characters of abelian differentials for some choice of 
complex structure $\tau=\tau(\chi)$ on $S$ consistent with the orientation. 
In other words, we describe the union
$$
\bigcup_{\tau\in T(S)} H^{1,0}(S_{\tau},\C),
$$
where $T(S)$ is the Teichm\"uller space of $S$. The proof is based upon 
Ratner's solution of Raghunathan's conjecture.  
\end{abstract}

\begin{dedication}
\hspace{4cm}
\vspace*{1cm}{To the memory of Sergei Kolyada}
\end{dedication}

\section{Introduction}

This paper is a slightly revised version of my preprint written in 2000 at the Max Plank Institute for Mathematics in Bonn. 
A few years after writing the preprint, I discovered a paper by Otto Haupt \cite{Haupt}, 
where the main result of my paper, Theorem 1.2 (including the genus 2 case), was proven by elementary methods. 
Another proof is contained in the  preprint of Bogomolov, Soloviev and Yotov, \cite{BSY} (who also study periods of pairs and even triples of abelian differentials). In view of Haupt's paper, the main point of my work is to establish a connection of the periods of abelian differentials to  ergodic theory. This connection and some of the methods used in this work were exploited 
by Calsamiglia,  Deroin and   Francaviglia in \cite{CDF} to further analyze the period map and to prove the connectivity of its fibers. In their paper they also found a mistake in my preprint, in the analysis of the genus 2 case, and gave a precise description of orbit closures in this setting. Therefore, I have removed the genus 2 case from the present paper;  otherwise, it remains essentially unchanged.

\medskip 
Let $S$ be a closed (i.e. compact with empty boundary)  
connected oriented surface of genus $n$. Recall 
that each complex structure $\tau$ on $S$ (consistent with the orientation) 
determines the linear subspace $H^{1,0}(S_{\tau},\C)\subset H^1(S,\C)$ of 
 complex dimension $n$ (i.e. half of the dimension of the cohomology 
group). In  down-to-earth terms, the subspace $H^{1,0}(S,\C)$ consists 
of the period characters of abelian differentials $\al\in \Om(S)$:
$$
\chi_\al=\chi \in H^1(S,\C), \quad \chi(c)= \int_c \al, \quad c\in H_1(S,\Z). 
$$
In this paper we describe the subset
$$
\bigcup_{\tau\in T(S)} H^{1,0}(S_{\tau},\C) \subset H^1(S, \C) 
$$
where $T(S)$ is the Teichm\"uller space of $S$. In other words,   
we give a necessary and sufficient condition for a character 
$\chi\in H^1(S,\C)$ to appear as the period of some abelian differential $\al$ 
on $S_\tau$ for {\em some} choice of the complex structure $\tau$ on $S$.

\begin{rem}
We note the difference between this question and the Schottky problem which 
asks for a description of the subvariety in the Grassmannian $G(n,2n)$ 
whose elements are subspaces  $H^{1,0}(S_{\tau},\C)$, with $\tau\in T(S)$. 
\end{rem}

Since the solution is obvious in the case $\chi=0$ we will consider 
only the nontrivial characters $\chi$. It turns out that there are 
precisely two topological obstructions for such $\chi$ to be the character 
of an abelian differential, the first is classical and is a part of 
the Riemann bilinear relations (see for instance \cite{Narasimhan}); 
the second is less known.

To describe the first obstruction (which applies for all $n\ge 1$) recall that 
the Poincar\'e duality defines a symplectic pairing 
$\om:H^1(S,\R)^{\otimes 2} \to \R$. This yields 
a quadratic form $H^1(S,\C)\to \R$ again denoted $\om$:
$$
\om(\chi):= \om(Re\chi, Im \chi).  
$$
If  $x_1, y_1,...,x_n,y_n$ denote the standard (symplectic) 
basis of $H^1(S,\Z)$ then $\om(\chi)$ equals
$$
\sum_{i=j}^n  Im(\ol{\chi(x_j)}\chi(y_j)). 
$$
The number $\om(\chi)$ can be also described as 
$$
\int_S f^*(dA),
$$  
where $dA$ is the area form $\frac{i}{2}dz \wedge \ol{dz}$ on 
$\C$, $f: S\to E$ is a section
of the complex line bundle $E$ over $S$ associated with $\chi$. 
(The form $dA$ is induced on $E$ via the projection 
$\tilde{S}\times \C \to \C$, where $\tilde{S}$ is the universal covering of 
$S$.)  

Note that in the case when $\chi\ne 0$ is the period character of an abelian 
differential $\al\in \Om(S)$ we have:
$$
\om(\chi)= \int_S \frac{i}{2}\al \wedge \bar\al
$$
is the area of the surface $S$ with respect to the singular Euclidean metric 
on $S$ induced by $\al$. Since this area has to be positive we get

\medskip
{\bf Obstruction 1.} If $\chi\in H^{1,0}(S_\tau)$ for some $\tau\in T(S)$ 
then $\om(\chi)>0$. 

\noindent The second obstruction applies only to special characters $\chi$ and surfaces of genus $n\ge 2$. 
In what follows we will regard elements of $H^1(S,\C)$ as additive 
characters $\chi$ on $H^1(S,\Z)$, this way we have the {\em image} of $\chi$, 
which is a subgroup $A_\chi$ of $\C$. 

\medskip
{\bf Obstruction 2.} Suppose that the image $\Image(\chi)$ of the character 
$\chi\in H^1(S,\C)$ is a discrete subgroup $A_\chi$ 
of $\C$ isomorphic to $\Z^2$ and $n\ge 2$. Thus $\chi$ gives rise to a homomorphism 
$$
\chi: H^1(S,\Z)\to H^1(T^2,\Z)
$$ 
where $T^2= \C/A_\chi$ is the 2-torus. This map is realized by a unique 
(up to homotopy) map $f: S\to T^2$. Then, for each 
$\chi\in H^{1,0}(S_\tau)$ the degree of $f$ has to be at least $2$. 

The reason for this obstruction is that if $\chi$ is the period of some 
$\al\in \Om(S_\tau)$ 
then the multivalued solution of the equation $dF=\al$ on the 
Riemann surface $S_\tau$ yields a (nonconstant) 
holomorphic map $f: S\to T^2$ which induces $\chi: H^1(S,\Z)\to H^1(T^2,\Z)$. 
Since the surface $S$ is assumed to have genus $n\ge 2$, 
the map $f$ cannot be a homeomorphism, hence its degree is at least $2$. 

Alternatively, the second obstruction can be described as follows. 
Assume again that the image $A_\chi$ of the character $\chi$ is a discrete 
subgroup isomorphic to $\Z^2$. Let $\Area(\chi)$ denote $\Area(\C/A_\chi)$, the 
area of the flat torus. Then the requirement $\deg(f)\ge 2$ is equivalent to
$$
\om(\chi) \ge 2 \Area(\chi). 
$$

We now assume that the surface $S$ has genus $n\ge 3$. Our main result is the following:

\begin{thm}\label{main}
If $n\ge 3$ and $\chi\in H^1(S,\C)$ satisfies the conditions imposed by the 1-st and 
the 2-nd obstruction then $\chi\in H^{1,0}(S_\tau)$ for some $\tau\in T(S)$. 
\end{thm}

\medskip 
In \S \ref{mero} we show that if $\chi$ is a nonzero character which is not 
the period of any abelian differential, it is nevertheless possible to 
find a complex structure $\tau$ on $S$ such that $\chi$ is the period 
character of a meromorphic differential with a single simple pole on $S_\tau$. 
We now identify the additive group $\C$ with the subgroup of $PSL(2,\C)$ 
consisting of translations $z\mapsto z+b, b\in \C$. Then we can regard $\chi$ as a representation 
$\rho: \pi_1(S)\to PSL(2,\C)$. For such $\rho$ define 
\begin{equation}
\label{for}
d(\rho):= {\bigg \{} \begin{array}{cc}
2n-2, &\hbox{~~if Obstructions 1 and 2 are satisfied,}\\
2n, & \hbox{~~otherwise.}
\end{array}
\end{equation}

We recall (see e.g. \cite{GKM}) that a {\em branched projective structure} $\si$ on a complex 
curve $S$ is an atlas with values in ${\mathbb S}^2$ where the local charts 
are nonconstant holomorphic functions (not necessarily locally 
univalent) and the transition maps are linear-fractional transformations 
(i.e. elements of $PSL(2,\C)$). Thus near each point $z\in S$ (which 
we identify with $0\in \C$) the local chart has the form $z\mapsto z^{m+1}$. 
The number $m=\deg(z)$ is called the degree of branching  at $z$. 
We get the {\em branching divisor} $D$ on $S$ whose degree is called the 
{\em degree of branching} $\deg(\si)$. For each representation 
$\rho: \pi_1(S)\to PSL(2,\C)$ there exists a complex-projective structure 
$\si$ (consistent with the orientation on $S$) which corresponds to {\em some}
 complex structure on $S$, such that $\rho$ is the holonomy of $\si$. 
We define $d(\rho)$ to be the least degree of branching for such structures. 
Note that for the trivial representation $\rho$, $d(\rho)=2n+2$ 
and the branched projective structure is given by the hyperelliptic covering. 
In this note we compute the function $d(\rho)$ in the very special 
case of representations with the image in the subgroup of translations. 
The general case will be treated elsewhere, here we only note that 
in  \cite{GKM} (see also \cite{Kapovich(1995)}) it was shown that for each 
representation $\rho$ with {\em nonelementary image}\footnote{I.e. the image 
does not have an invariant finite nonempty subset in 
$\H^3\cup {\mathbb S}^2$.}, $d(\rho)\in \{0,1\}$ equals the 2-nd 
Stiefel--Whitney class of $\rho$ (mod 2).

\begin{cor}\label{corfor}
Suppose that $n\ge 3$. For each nontrivial 
representation $\rho: \pi_1(S)\to PSL(2,\C)$ whose image is contained in 
the subgroup of translations, the function $d(\rho)$ is given by the formula 
(\ref{for}).  
\end{cor}

The lower bounds in this theorem are given by the Riemann--Roch theorem (see \S 
\ref{mero}), while the upper bound follows from Theorems \ref{main} 
and \ref{pole}.

Since the map $P: \al\to \chi_\al$, which sends the abelian differential to its 
character, is complex-linear, it suffices to prove Theorem \ref{main} for 
{\em normalized} characters, i.e. the characters $\chi$ such that 
$\om(\chi)=1$ (hence the 1-st obstruction automatically holds). 
We let
$$
X:= \{ \chi\in H^1(S,\C) : \om(\chi)=1\}
$$ 
and 
$$
\Si:= X\cap \bigcup_{\tau\in T(S)} H^{1,0}(S_{\tau},\C). 
$$
Let $\Om$ denote the vector bundle over $T(S)$ whose fiber over a 
point $\tau\in T(S)$ consists of abelian differentials $\Om(S_\tau)$. 
We let $\Om'$ denote the submanifold in $\Om$ consisting of abelian 
differentials $\al$ such that $\om(\al)=1$.   
We have the map 
$$
P: \Om' \to \Si\subset X. 
$$
To explain the appearance of  ergodic theory in the proof we will need 
two elementary facts about the subset $\Si$ in $X$. 

{\bf Fact 1.} (See \S \ref{TH}.) 
The map $P: \Om'\to X$ is open. In particular,  $\Si$ is open in $X$.

\medskip
We let $G=Sp(n)=Sp(2n,\R)$ denote the group of linear symplectic automorphisms 
 of the symplectic structure $\om$ on $\R^{4n}= H^1(S,\C)$. This 
is a simple algebraic Lie group which acts naturally on $X$. 
It is elementary that the action of $G$ on $X$ is transitive. 
The stabilizer $G_\chi$ of a point $\chi\in X$ is isomorphic to 
$Sp(2n-2)$. Thus $X=Sp(2n)/Sp(2n-2)$. Recall that the integer 
symplectic group $\Ga=Sp(2n,\Z)$ is a {\em lattice} in the group $G$. 

\medskip
{\bf Fact 2.} The subset $\Si$ is invariant under $\Ga$. 

Recall that the group of 
orientation-preserving diffeomorphisms $\Diff(S)$ acts on $H^1(S,\C)$ 
through the group $\Ga$. If $\chi\in \Si$ is the period character 
of $\al\in \Om(S_\tau)$ and $\ga\in\Ga$ corresponds to a diffeomorphism 
$h:S\to S$, then $\ga(\chi)$ is the period character of the abelian 
differential 
$$
h^*(\al)\in \Om(S_{h^*(\tau)}),
$$ 
where $h^*(\tau)$ is the pull-back of the complex structure $\tau$ via $h$. 
Thus $\ga(\Si)= \Si$. 

Combining the above two facts we see that $\Si$ is a (nonempty) open 
$\Ga$-invariant subset of $X$. We recall 

\begin{thm}
[C. Moore, see \cite{Zimmer(1984)}]  
If $G$ is a semisimple Lie group, $\Ga$ is a lattice in $G$ and 
$H$ is a noncompact Lie subgroup in $G$ then $H$ acts ergodically on 
$\Ga\backslash G$. Equivalently, $\Ga$ acts ergodically on $G/H$. 
\end{thm}

Thus, since $\Si\subset X=Sp(2n)/Sp(2n-2)$ is an open nonempty 
$\Ga$-invariant subset, the complement 
$X-\Si$ has zero measure. In particular, $\Si$ is dense in $X$. 
Ergodicity of the action $\Ga\acts X$ implies that 
{\em generic}\footnote{In the 
measure-theoretic sense.} points $\chi\in X$ have dense $\Ga$-orbits.  
Our objective is to understand the {\em nongeneric} orbits. 
This is done by applying Ratner's solution of Raghunathan's conjecture. 
Ratner's theorem implies that there are only few types of nongeneric orbits. 
We will show that most of them correspond to the characters with discrete 
image. After we describe other orbits we will show that Obstruction 2 suffices 
for the existence of an abelian differential with the given period character.

\medskip
{\bf Acknowledgments.} During the period of this work I was partially 
supported by the National Science Foundation 
and by the Max--Plank Institute (Bonn) for whose hospitality I am grateful. I am also grateful to Bertrand Deroin for pointing out 
 my mistake in the genus two case and to Curt McMullen for repeatedly asking me to publish this work. I am thankful to Pieter Moree and to the anonymous referee for numerous suggestions and corrections.

\section{Geometric preliminaries}
\label{TH}

{\em Geometric interpretation of nonzero abelian differentials $\al$.} 
Each nonzero abelian differential $\al\in \Om(S_\tau)$ 
determines a singular Euclidean structure on the surface $S$ 
with isolated singularities at zeroes of $\al$, see \cite{Strebel}. 
Let $\Zero(\al)\subset S$ denote the set of zeroes of $\al$. 

The local charts for this structure are given by the branches of the 
indefinite integral
$$
F(z)=\int_{z_0}^z \al,
$$
where $z_0\in S$ is a base-point. 
If $\al$ vanishes (at the order $m-1$) at a point $0\in S$ then the 
local chart at $0$ is a $k$-fold ramified covering $z\mapsto z^{m}$. 
The transition maps of the flat atlas on $S- \Zero(\al)$ are Euclidean 
translations. Vice-versa, suppose that we are given a flat structure 
on the (topological) surface $S$ where the local charts have the form  
$z\mapsto z^{m}$, $m\ge 1$, and the transition maps away from the 
branch-points are Euclidean translations. This structure canonically 
defines a complex structure on $S$ together with an abelian differential 
$\al$ obtained by the pull-back of $dz$ via the local charts. Every such 
singular Euclidean structure gives rise to a {\em developing map} 
$dev: \t S\to \C$ where $\t S$ is the universal abelian covering of $S$ and 
$H:=H_1(S,\Z)$ acts on $\t S$ by deck-transformations. The mapping $dev$ is 
$\chi$-equivariant, where $\chi: H_1(S,\Z)\to \C$ is the holonomy 
of the above structure (it coincides with the character of the 
associated abelian differential). 
The space $E(S)$ of the above Euclidean structures has a natural topology: the 
topology of  uniform convergence on compacts 
of the developing mappings. It is easy to see 
that with this topology the natural bijection $E(S)\to \Om- 0_\Om$ is 
a homeomorphism.   Here $0_\Om$ is the image of the zero-section of the bundle $\Om\to T(S)$, i.e. 
$0_\Om$ consists of zero abelian differentials.

\medskip
{\bf Matrix form of the characters.} Given the standard (symplectic) 
basis in $H_1(S,\Z)$, $x_1,y_1,...,x_n,y_n$, we can identify each 
character $\chi: H_1(S,\Z)\to \C=\R^2$ with the $2\times 2n$ matrix 
$$
M(\chi):= [M_1 M_2 ... M_n], 
$$
$$
M_j=M_j(\chi):=  \left[ \begin{array}{cc}
a_j & b_j\\c_j & d_j \end{array}\right], \quad j=1,...,n.
$$
Here 
$$
\chi(x_1,...,y_n)= (u,v)^t, \quad u=(a_1,b_1,...,a_n,b_n), \quad v=(c_1,d_1,...,c_n,d_n),  
$$
and the vectors $u,v$ are the row-vectors of the matrix $M$. The 
group $G=Sp(2n)$ acts on the matrices $M$ by multiplying them from the right.  
The matrix $M(\chi)$ is the {\em matrix form} of the character $\chi$. 
Then we define
$$
\om_j(u,v)=det(M_j(\chi))= 
\left|\begin{array}{cc}
a_j & b_j\\
c_j & d_j
\end{array}\right|, \ j=1,...,n;
$$
it follows that $\om(u,v)=\sum_j \om_j(u,v)$. The group $SL(2)=Sp(2)$ acts 
on the characters $\chi$ by multiplying their matrices from the left. 
It is clear that this action commutes with the action of 
$Sp(2n,\Z)\subset G$ and that it preserves each determinant $\om_j(\chi)$.

\begin{lem}
\label{rem1}
$Sp(2)\Si=\Si$. 
\end{lem}
\proof Suppose that $\chi\in\Si$ is the period character of an abelian 
differential corresponding to a singular Euclidean structure $\si$. 
Take $A\in Sp(2)$. Composing coordinate charts of $\si$ with $A$ deforms 
$\si$ to a new singular Euclidean structure of the same area. 
The holonomy of this structure is the composition $A\circ \chi$. 
Hence $A\chi\in \Si$. \qed

\begin{lem}
\label{rat}
Suppose that $\chi=(u,v)$ and $u,v\in \R^{2n}$ span a 
2-dimensional rational subspace (i.e. a subspace which admits a rational 
basis). Then the $\Z$-module ${\mathcal M}$ generated by the columns of the 
matrix $M(\chi)$ has rank 2, i.e. is discrete as a subgroup of $\R^2$.     
\end{lem} 
\proof The action of $GL(2)$ by  multiplication from the left on the 
matrix $M(\chi)$ preserves the rank of ${\mathcal M}$. Since 
$\Span(u,v)$ is a rational subspace there exists 
a matrix $A\in GL(2)$ such that the matrix $AM(\chi)$ has integer entries. 
The rank of the $\Z$-module generated by its columns is clearly $2$. 
\qed

\medskip
Define 
$$
X_+:= \{ \chi\in X: \om_j(\chi)>0, j=1,...,n\}. 
$$
Our strategy in dealing with the {\em nongeneric characters} 
$\chi\in X$ is to find $\ga\in Sp(2n,\Z)$ such that $\om_j(\ga\chi)>0$,  
 $j=1,...,n$, i.e. $\ga\chi\in X_+$. 
As we will see in Theorem \ref{geo} the existence of such $\ga$ 
would imply that $\chi$ belongs $\Si$  (i.e. that $\chi$ is the period 
character of an abelian differential).

\begin{figure}[tbh]
\centerline{\epsfxsize=5in \epsfbox{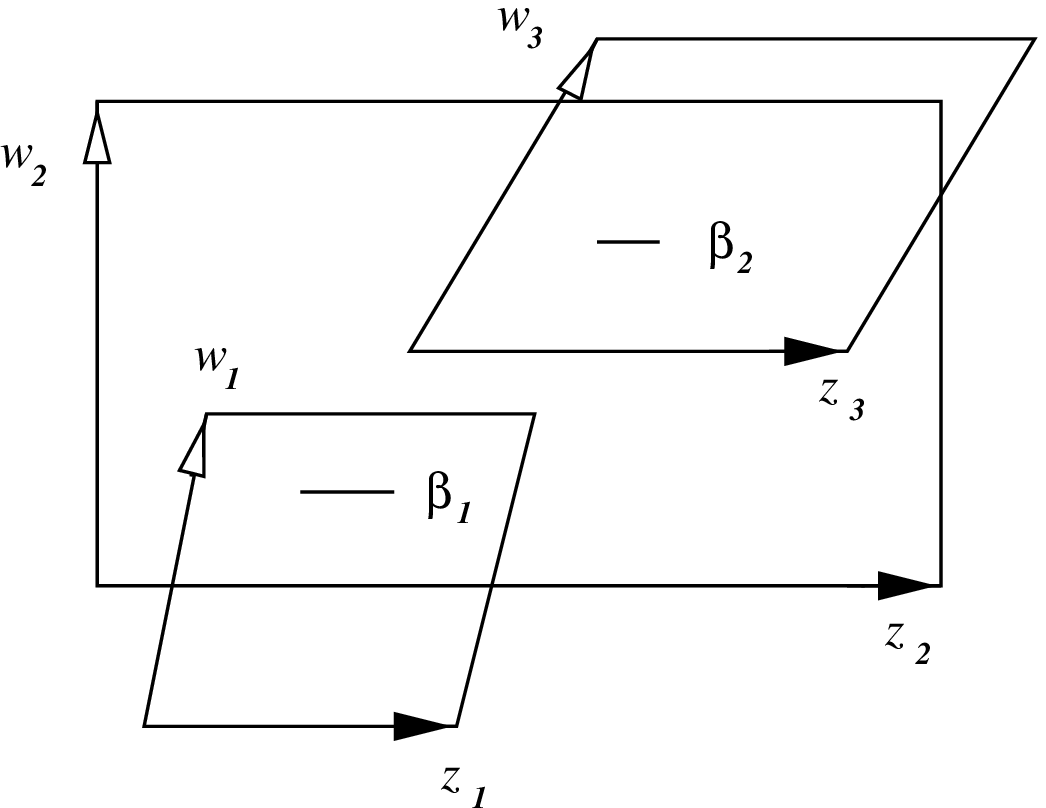}}
\caption{\sl }
\label{F}
\end{figure}

\begin{thm}\label{geo}
$X_+\subset \Si$. 
\end{thm}
\proof Let $(u,v)\in X_+$, $u=(a_1,b_1...,a_{n},b_{2n}), 
v=(c_1,d_1,...,c_n,d_n)$. 
We let $z_j:=(a_j,c_j), w_j:=(b_j,d_j)\in \R^2$, $j=1,...,n$. 
Each pair of vectors $(z_j,w_j)$ determines a fundamental parallelogram $P_j$ 
in $\R^2$ for the lattice generated by $z_j,w_j$. 
Using parallel translations place these parallelograms 
such that $P_j\cap P_{j+1}$ has nonempty interior, $j=1,...,n-1$. Then for each 
pair of parallelograms $P_j, P_{j+1}$ ($j=1,...,n-1$) 
cut both $P_j, P_{j+1}$ open along common segments $\be_j$   and then 
glue them along the resulting circles. Call the result $\Phi$. 
See Figure \ref{F}.

Finally, for each parallelogram 
$P_j$ identify the opposite sides via a parallel translation. The result is 
a surface $S$, equipped with the projection $\del: \tilde{S}\to \C$ where 
$\tilde{S}$ is the universal abelian covering. The surface $\Phi$ is the 
fundamental domain for the action of $H_1(S,\Z)$ on $\t S$ via deck 
transformations. The restriction $\del|\Phi : \Phi\to \C$ is the 
obvious projection. Note that $\del$ is a local homeomorphism away from 
the translates of the end-points of the segments 
$\be_j$. Near the end-points of such segments the mapping $\del$ is a 2-fold 
ramified covering. The abelian differential $\al$ on 
$S$ is obtained by taking the pull-back of $dz$ from $\C$ to $\t S$ via $\del$ 
and then projecting it to $S$. The edges 
of the parallelograms $P_j$ correspond to the standard generators of 
$H_1(S,\Z)$. It is clear that the periods of $\al$  over the generators 
of $H_1(S,\Z)$ are given by evaluation of $\chi$ on these generators. \qed 

\medskip
The above lemma implies that it suffices to show that 
$\Ga\chi\cap X_+\ne \emptyset$ to prove that $\chi\in \Si$. 
Note however that there are characters in $\Si$ which do not belong to 
the orbit $\Ga X_+$. These are the characters with the discrete image 
$A_\chi\cong \Z^2$ such that 
$$
\frac{\om(\chi)}{\Area(\C/A_\chi)}< n. 
$$

To find abelian differentials corresponding to such characters we need 
another construction that we describe below.

\begin{figure}[tbh]
\centerline{\epsfxsize=5in \epsfbox{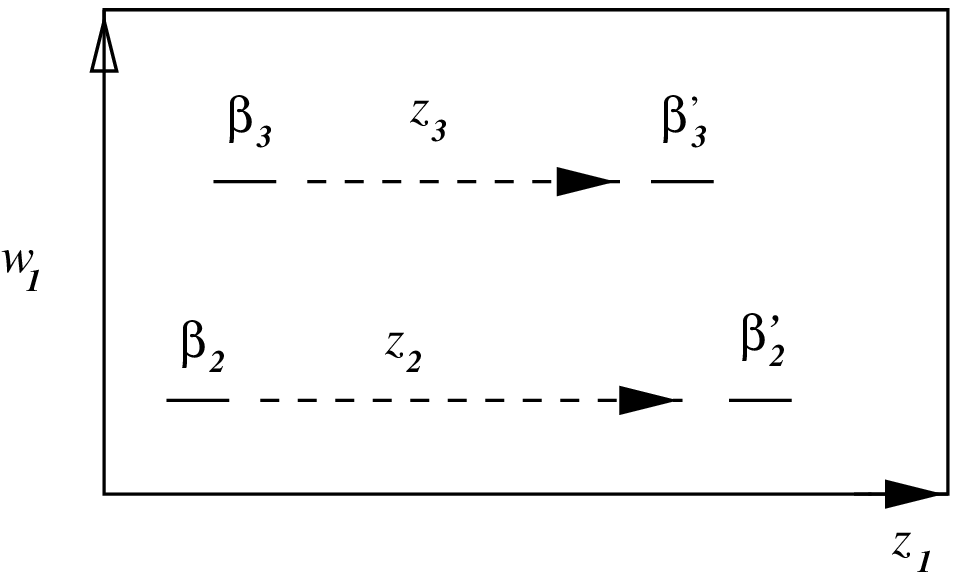}}
\caption{\sl }
\label{F2}
\end{figure}

\begin{lem}
Suppose that the character $\chi$ has the matrix form
$$
[M_1 M_2 ... M_n], 
M_1= \left [\begin{array}{cc}
a_1=\om(\chi) & 0\\0 & 1 \end{array}\right], 
M_j= \left[\begin{array}{cc}
a_j & 0\\0 & 0 \end{array}\right], j=2,...,n, 
$$ 
where $0< a_j< a_1$, $j=2,...,n$. 
Then $\chi\in \Si$. 
\end{lem}
\proof Similarly to the previous lemma we construct a complex structure 
and an abelian differential by gluing certain polygons. Let $P_1$ be the 
fundamental rectangle for the group generated by the vectors $z_1,w_1$ which 
are the columns of $M_1$. Inside $P_1$ choose pairwise 
disjoint horizontal segments $\be_j,\be_j', j=2,..,n,$ such that the translation via 
$[a_j 0]$ sends $\be_j$ to $\be_j'$. We then cut $P_1$ open along the segments 
$\be_j,\be_j'$ and identify the resulting circles via the translations by 
$[a_j 0]$, $j=2,..,n$. Finally, glue the sides of $P_1$ via the 
horizontal translations, see Figure \ref{F2}. 
Analogously to the previous lemma we get a 
singular Euclidean structure with the holonomy $\chi$. The singular 
points of this structure correspond to the end-points of the segments 
$\be_j$ (the total angle at each of these points is $4\pi$).   
\qed

\begin{lem}\label{completion}
Suppose that $u,v\in \Z^4$ are vectors such that $\om(u,v)=1$. Then 
this pair of vectors can be completed to an integer symplectic basis 
in $\R^4$. 
\end{lem}
\proof Let $W:= \Span(u,v)$. Recall that the symplectic projection 
$\Proj_W(z)$ of
a vector $z$ to $W$ is given by
$$
\Proj_W(z)= \om(z,v)u- \om(z,u)v.
$$
Hence $\ker(\Proj_W)=W^\perp$ is a rational subspace in $\R^4$ and we choose a basis 
$p,q\in W^\perp$ such that the vectors $p,q$ generate the abelian group 
$\Z^4\cap W^\perp$. The vectors $u,v,p,q$ generate the group $\Z^4$ since the 
symplectic projection of $\Z^4$ to $W$ and $W^\perp$ is contained in 
$\Z^4\cap W$ and $\Z^4\cap W^\perp$ respectively. It follows that $\om(p,q)=1$ and 
$x,y,p,q$ form an integer symplectic basis in $\R^4$. \qed

\begin{lem}\label{zero}
Suppose that $u\in \R^{2n}$ is a nonzero vector. 
Then there exists $\ga\in \Ga$ such that no coordinate of $\ga(u)$ 
is zero. If $u,v\in \R^{2n}$  are such that $\om(u,v)>0$, then 
there exists $\ga\in \Ga$ such that $\om_j(\ga(u),\ga(v))\ne 0$ for 
each $j=1,...,n$. 
\end{lem}
\proof The projection $Sp(2n)\to \R^{2n}-0$ given by $g\mapsto g(\ov{e_1})$ 
is a real algebraic morphism. The union 
$$
\bigcup_{j=1}^n \{ x\in \R^{2n} : x_j=0\}
$$
is a proper (real) algebraic subvariety, hence its inverse image $Y$ 
in $G=Sp(2n,\R)$ is again a proper algebraic subvariety. Since $\Ga$ is 
Zariski dense in $G$ we conclude that $Y$ is not $\Ga$-invariant. 
The proof of the second assertion is similar and is left to the reader. 
\qed

Recall that $\Om$ denotes the vector bundle over the Teichm\"uller space  
$T(S)$ where the fiber over a point $\tau$ consists of abelian differentials 
on the Riemann surface $S_\tau$; $0_\om$ denotes the image of the zero section 
of $\Om$.  We have the period map $P: \Om\to H^1(S,\C)$, 
$\al\mapsto \chi_\al$. 

The following theorem is a variation on the Hejhal--Thurston Holonomy theorem, 
see \cite{Hejhal}, \cite{Thurston(1978-81)},  and 
\cite{Epstein-Canary-Green}, \cite{Goldman(1987b)}. See also 
\cite[Section 12]{GKM} for an alternative argument. 

\begin{thm}
\label{holonomy map}
(The Holonomy Theorem.)  
The restriction mapping $P: \Om - 0_\Om\to H^1(S,\C)$ is open. 
\end{thm}
\proof To prove this theorem we be using  a geometric description of 
the nonzero abelian differentials $\al$ given in the beginning of this section.    
Let $\si\in E(S)$ be a singular Euclidean 
structure with the period character $\chi$. Let $f: \t S\to \C$ denote 
the developing mapping of $\si$. Suppose that 
$\chi_k: H_1(S,\Z)\to \C$ is a sequence of characters converging to $\chi$. 
Our goal is to find (for large $k$) structures $\sigma_k\in E(S)$ with the  
period characters  $\chi_n$ and such that $\lim_k \sigma_k =\sigma$. 

Choose a triangulation $T$ of $S$ such that each edge is a geodesic arc 
with respect to the singular Euclidean structure $\si$ and each simplex is 
contained in a coordinate neighborhood of $\si$. We will assume that each 
singular point of $\si$ is a vertex of this triangulation. 
Lift this triangulation to a triangulation $\t{T}$ of $\t{S}$ of $S$. 
Pick a finite collection $\Delta_1,...,\Delta_M$ 
of 2-simplices in $\t{T}$, one for each $H$-orbit. 
Let $g_i, i=1,...,N$, be the elements of 
the deck-transformation group $H$, such that 
$$
g_i(\cup_j \Delta_j)\cap 
\cup_j \Delta_j\ne \emptyset.
$$
 Let $C$ be a compact subset of $\t{S}$ 
whose interior contains  both $D:=\cup_j \Delta_j$ and its images under 
$g_i$'s. For each $\chi_k$ we construct a continuous 
$\chi_k$-equivariant mapping $f_k: D\to \C$ such that: 

(i) $f_k$ maps each 2-simplex homeomorphically to a Euclidean 
2-simplex in $\C$.

(ii) $f_k$'s  converge to $f|D$ uniformly on compacts. 

Finally, extend each $f_k$ to a $\chi_k$-equivariant mapping  
$f_k: \t S\to \C$. It remains to show that each mapping 
$f_k$ is a local homeomorphism for large $k$ (away from the singular points) 
and is the $m(x)$-fold ramified covering at each point where $f$ is such 
a covering. It suffices to check this for points in $D$. 

(a) If $x\in int(C)$ belongs to the 
interior of a 2-simplex in $\cup_i g_i D$, then the claim follows since each 
$f_n$ is a homeomorphism on each simplex. 

(b) Suppose $x$ belongs to the 
interior of a 
common arc $\eta$ of two 2-simplices $\Delta,\Delta'$ in $\cup_i g_i D$. 
Since $f$ is a local homeomorphism, $f(\Delta), f(\Delta')$  lie (locally) 
on different sides of the  segment $f(\eta)\subset \C$. Therefore the 
same holds for $f_k$ if $k$ is sufficiently large. Thus, $f_k$ 
does not ``fold'' along the arc $\eta$ 
and is a local homeomorphism at $x$.  

(c) Lastly, if $x$ is a vertex of a simplex, then 
the degree of $f$ at $x$ equals $m(x)$, hence for large $k$, 
the degree of $f_k$ at $x$ is $m(x)$ and it follows from (b) that 
$f_k$ is a $m(x)$-fold ramified covering at $x$. 

Equivariance of $f_k$'s implies that they converge to $f$ 
uniformly on compacts. \qed 

\medskip
{\bf Line stabilizers in $Sp(2n)$}. In what follows we will need a 
description of the subgroups $B$ in $Sp(2n)$ with invariant line 
$L\subset \R^{2n}$. Let $V\subset \R^{2n}$ be a 2-dimensional 
symplectic subspace containing $L$. To describe the 
structure of the group $B$ we have to recall several facts about  
{\em Heisenberg groups}. Consider the $2n-2$-dimensional 
symplectic vector space $(V, \om|V)$. The Heisenberg group 
corresponding to this data is the $2n-1$-dimensional Lie group 
which fits into the short exact sequence
$$
1 \to \R \to H_{2n-1} \to V\to 1, 
$$
where $V$ is treated as the abelian (additive) Lie group. The normal subgroup 
$\R$ is central in $H_{2n-1}$. If $g,h\in H_{2n-1}$ project to the vectors 
$x,y\in V$ then $[g,h]=\om(x,y)\in \R$. The {\em Heisenberg dilation} 
on this group is the action of the (multiplicative) group  $\R_+$ on 
$H_{2n-1}$ such that $t\in \R_+$ acts on the center $\R\subset H_{2n-1}$ 
via multiplication by $t^2$ and acts on $V$ via multiplication by $t$. 
Given this one defines the Lie group $H_{2n-1}\rtimes \R_+$ where 
$\R_+$ acts on the Heisenberg group via Heisenberg dilation. 
One can show that the resulting Lie  group acts simply-transitively 
on the complex-hyperbolic space $\C\H^{n}$ of the complex dimension $n$, 
however we will not need this fact. What we will use is the following 
elementary lemma. 

\begin{lem}\label{flow}
The $2n$-dimensional Lie group $CH_{2n}:=H_{2n-1}\rtimes \R_+$ 
contains no lattices. 
\end{lem} 
\proof Suppose that $\Del$ is a discrete subgroup of 
$H_{2n-1}\rtimes \R_+$ with the quotient $M=H_{2n-1}\rtimes \R_+/\Del$. 
The unit speed flow on $H_{2n-1}\rtimes \R_+$ 
along the $\R_+$-factor is volume-expanding and $\Del$-invariant. 
Hence it yields a volume-expanding flow on $M$. It follows that 
$\Vol(M)=\infty$. \qed

\medskip 
We are now ready to describe the structure of $B$. The group $B$ preserves 
the span $L+V$ of $L$ and $V$, the projection $L+V\to V$ along the $L$-factor 
transfers the action of $B$ to the action of the symplectic group $Sp(2n-2)$ 
on $V$. The kernel of the homomorphism $B\to Sp(2n-2)$ is the group 
$CH_{2n}=H_{2n-1}\rtimes \R_+$. Here the $\R_+$-factor acts trivially 
on $V$ and as the maximal torus in $Sp(2)\acts V^\perp$ preserving $L$. 
The center $\R$ of the Heisenberg group $H_{2n-1}$ is the kernel of the 
action $B\acts L+V$. The whole group $B$ splits as the semidirect product 
$CH_{2n}\rtimes Sp(2n-2)$, where $Sp(2n-2)$ acts by conjugation on the 
$V$-factor of $H_{2n-1}$ the same way it acts on the vector space $V$. 
The subgroup $Sp(2n-2)$ commutes with the subgroup $B_0:=\R\rtimes \R_+$, 
where $\R$ is the center of $H_{2n-1}$. The proof of these assertions is a 
straightforward  linear algebra computation and is left to the reader. 

\begin{defn}
The group $H_{2n-1}$ is called the {\em Heisenberg group} associated to the 
flag $(V,L)$ in $(\R^{2n},\om)$, where $V$ is a 2-dimensional 
symplectic subspace and $L$ is a line. 
\end{defn}

\section{Ratner's Theorem}

Let $G$ be a reductive algebraic Lie group and  
$U\subset G$ be a connected subgroup generated by unipotent 
elements\footnote{I.e. elements 
whose adjoint action on the Lie algebra of $G$ is unipotent.}.  
Suppose $\Ga\subset G$ is a lattice, i.e. a 
discrete subgroup with the quotient $\Ga\backslash G$ of finite volume 
(with respect to the left-invariant measure on $G$). Important examples of 
lattices in algebraic Lie groups $G$ defined over $\Q$ are given by the 
{\em arithmetic groups}, i.e. subgroups commensurable with $G_\Z$, 
the group of integer points in $G$. The group $U$ acts by right 
multiplications on the manifold $M=\Ga\backslash G$. On the other hand, 
the group $\Ga$ acts by the left multiplication on the manifold $X=G/U$. 
Given $g\in G$ we let $[g]$ denote its projection to $M$. 

\begin{thm}
[M. Ratner, see \cite{Ratner(1991),Ratner(1995)}]  
Under the above conditions for each $g\in G$ the closure (in the 
classical topology) of $[g]U$ in $M$ is ``algebraic''. 
More precisely, there exists a Lie subgroup $H\subset G$ such that 

\begin{itemize}
\item $\ol{[g]U}= [g]H$. 

\item $H^g\cap \Ga$ is a lattice in $H^g:= gH g^{-1}$. 
\end{itemize} 
\end{thm}

This result is known as {\em Raghunathan's Conjecture}. 
Special cases of this conjecture were proven before Ratner by Dani 
\cite{Dani(1986)} and Margulis \cite{Margulis(1986)}. Actually, Ratner's 
theorem does more than what is stated above: it describes 
$\Ga$-invariant ergodic measures on $M$ and uses 
the ergodic framework to prove Raghunathan's Conjecture. We note that 
the group $H$ may not be connected, however if $H(0)$ is the connected 
component of the identity in $H$, then $H(0)\cap \Ga$ is still a lattice 
in  $H(0)$. Below we reformulate Ratner's theorem in terms of the action 
of $\Ga$ on $G/U$. Let $g\in G$ be the element which projects to $x$. Then
$$
\ol{\Ga g U}= \Ga g H= \Ga H^g g.
$$
Hence 

\begin{cor}
Suppose that $X:= G/U$ and $x=gU\in X$. Then the closure of $\Ga x$ in $X$ 
equals the $H^g$-orbit of $x$ in $X$, where $H^g$ is a Lie subgroup of $G$ 
such that $H^g\cap \Ga$ is a lattice in $H$. 
\end{cor}

Note that $gUg^{-1}=G_x$ is the stabilizer of $x$ in $G$. 
By taking the connected component of the identity we get: 

\begin{cor}
The closure $\ol{\Ga x}$ in $X$ contains the orbit $\Ga F_x x$, where 
$F_x$ is a connected Lie subgroup of $G$ which contains $G_x$ 
and $\Ga\cap F_x$ is a lattice in $F_x$. 
\end{cor}

Ratner's theorem gives a tool for describing the {\em exceptional} orbits 
for the $\Ga$-action on $X$, still, some work has to be done by 
analyzing various Lie subgroups $F_x\subset G$ which might appear.

We now specialize to the case $G= Sp(2n, \R)$, the automorphism group 
of the standard symplectic form $\om$:
$$
\om(a_1, b_1,...,a_{n},b_{n})= \sum_{j=1}^{n} a_j b_{j+1} - a_{j+1}b_j,  
$$
and $X\subset (\R^{2n})^2$ consists of the pairs of vectors $u,v$ such 
that $\om(u,v)=1$. 

The stabilizer $U$ of the point $(\ov{e_1},\ov{e_2})\in X$ is the group 
$Sp(2n-2,\R)$ embedded in $G$ as the subgroup of block-diagonal matrices:
$$
\left[\begin{array}{ccc}
1 & 0 & 0  \ldots 0\\
0 & 1 & 0  \ldots 0\\
0 & 0 & Sp(2n-2)
\end{array}\right].
$$
Although the group $U$ is not unipotent itself, it is generated by 
 unipotent elements, hence Ratner's theorem applies. 
Recall that $\Ga=Sp(2n,\Z)$ is a lattice in $G$, we also note that 
$\Ga \cap U$ is a lattice in $U$ as well. 
In the rest of the paper we will use the notation  
$U'=G_\chi$ to denote the stabilizer of the point $\chi\in X$.  

\medskip
{\bf Connected Lie subgroups of $G$ containing $U$.} 
To apply Ratner's theorem we have to know which Lie subgroups of $G$ 
contain the Lie subgroup $U'$ (conjugate to $U$).  We will 
list all {\em maximal} subgroups containing $U$. 
Recall that a connected Lie subgroup $G_1\subset G$ 
is said to be {\em maximal} if it is not 
contained in any proper connected 
Lie subgroup $G_2\subset G$. We will use a classification of 
maximal subgroups of classical complex Lie groups done by Dynkin 
\cite{Dynkin(1952)} (the real case was carried out by Karpelevich 
\cite{Karpelevich(1955)}). In our case the classification 
of maximal subgroups of $Sp(2n,\C)$ easily implies (via the complexification) 
the needed result for the group of real points $Sp(2n,\R)$.

\begin{thm}
[E. Dynkin, see Ch. 6, Theorems 3.1, 3.2 in \cite{OV}]  
Suppose that $H\subset Sp(2n,\C)$ is a maximal connected Lie 
subgroup. Then one of the following holds:

(a) $H$ is a maximal parabolic subgroup of $Sp(2n,\C)$. 

(b) $H$ is conjugate to the subgroup $Sp(k,\C)\times Sp(N-k,\C)$. 

(c) $H$ is conjugate to $Sp(s,\C)\otimes SO(t,\C)$ where $2n=st$, 
$s\ge 2, t\ge 3$, $t\ne 4$ or $s=2,t=4$. 
\end{thm}

Note that in our situation $H$ contains $U\cong Sp(2n-2,\C)$, 
hence we can ignore the case (c). In the case (b) the only 
possibility is that $F$ is conjugate 
to the group $Sp(2,\C)\times Sp(2n-2,\C)$. In the case (a) the group $H$ 
has to preserve a complex line in $\C^{2n}$. 

We let $\chi= (u,v)$, $u,v\in \R^{2n}$ 
are such that $\om(u,v)=1$. Let $V$ denote $\Span(u,v)$. 
The group $U'=G_\chi\cong Sp(2n-2,\R)$ fixes the vectors $u,v$. 
This group also acts as the full group of linear symplectic automorphisms 
of the symplectic complement $V^\perp\cong \R^{2n-2}$ of $V$.  
The maximal subgroups of $G$ which contain $U'$ are:

\begin{enumerate}
\item The group $H=Sp(V)\times U'$, where $Sp(V)\cong Sp(2,\R)$ is the group 
of automorphisms of $V$. (The semisimple case.) 

\item The maximal parabolic subgroup $H$ of $G$ which has an invariant line 
$L\subset \R^{2n}$. (The non-semisimple case.) We note that in this 
case $L$ is necessarily contained in $V$. 
\end{enumerate}

Recall that in each case we have to find connected subgroups $F_\chi\subset 
H$ which contain $G_\chi=U'$ and such that $F_\chi\cap \Ga$ is a 
lattice in $F_\chi$.

\section{The semisimple case} 

In this case the group $F_\chi \subset Sp(V)\times U'$ 
containing $U'$, splits as the direct product 
$$
F_\chi\cong S \times Sp(2n-2,\R), 
$$
where $S\subset Sp(2,\R)$. We will need the following

\begin{theorem}
[See e.g. \cite{Margulis(1991)}] 
Suppose that $F_1,F_2$ are simple real algebraic Lie groups 
such that their complexifications do not have isomorphic Lie algebras. 
Then any lattice $\Del\subset F_1\times F_2$ is reducible, i.e. 
$\Del\cap F_i$ is a lattice for each $i=1,2$. 
\end{theorem}

We also recall (see \cite[Corollary 8.28]{Raghunathan}): 

\begin{thm} 
[M. Raghunathan, J. Wolf]  \label{rw}
Suppose that $F$ is a connected Lie group whose semisimple part contains no 
compact factors acting trivially on the radical $R(F)$ of $F$. Then 
each lattice $\Del\subset F$ intersects the radical $R(F)$ along a sublattice 
in $R(F)$. Moreover, the projection of $\Del$ to $F/R(F)$ is a lattice 
in this Lie group. 
\end{thm}

In our case the group $S$ is either solvable or equals $Sp(2)$, hence 
combining the two above theorems we conclude that either:

(i) $\Ga\cap U'$ is a lattice, or 

(ii)  $n=2$, $F_\chi\cong Sp(2)\times Sp(2)$ and  $\Ga\cap U'$ is not a 
lattice\footnote{We note that the group $Sp(2)\times Sp(2)$ 
contains irreducible lattices, namely the Hilbert modular groups.}. 

In view of the assumption that $S$ has genus $\ge 3$, we are considering here only case (i),  
when $\Ga\cap U'$ is a lattice.

By the Borel density theorem (see e.g. \cite{Zimmer(1984)}) the 
intersection $U'\cap\Ga$ is Zariski dense in $U'$, in particular 
it contains a diagonalizable matrix 
$A\in Sp(2n)$ which has the eigenvalue $1$ of the multiplicity $2$. 
Since $A$ has rational 
entries, the kernel $\ker(A-I)$ is a rational subspace. 
We recall that the group $U'$ is the pointwise stabilizer of the linear 
subspace $\Span(u,v)$ of $\R^{2n}$ spanned by $u=Re(\chi), v=Im(\chi)$. 
Hence $\Span(u,v)$ is a rational subspace of $\R^{2n}$. 

Lemma \ref{rat} thus implies that the image $A_{\chi}$ of the character 
$\chi: H_1(S,\Z)\to \C$ is a discrete subgroup of $\C$ isomorphic to $\Z^2$. 
Moreover, without loss of generality we can assume that $A_{\chi}$ 
is the standard integer lattice in $\C$ (see Section \ref{TH}). 
This might require scaling $\om(\chi)$ by a positive real number.

We recall that $\om(u,v)>0$, where $\chi=(u,v)$,
$$
u=(a_1,b_1,...,a_n,b_n), \quad v=(c_1,d_1,...,c_n,d_n), \quad a_j,b_j,c_j,d_j\in \Z. 
$$ 

\begin{lem}
There exists $\ga\in \Ga$ such that the character 
$\ga\chi=\chi'=(u',v')$ satisfies:

(i) $\om_1(u',v')>0$.

(ii) $\om_j(u',v')=0$ for each $j\ge 2$ and, moreover, 
$$
M_j(\chi')= \left[\begin{array}{cc}
a_j' & b_j'\\
c_j' & d_j'
\end{array}\right]= 
\left[\begin{array}{cc}
a_j' & 0\\
0 & 0
\end{array}\right], a_j' \ge 0. 
$$
\end{lem}
\proof We recall that without loss of generality we can start with $(u,v)$ 
such that for each $j=1,...,n$, $\om_j(u,v)\ne 0$ or 
$$
M_j(\chi)= \left[\begin{array}{cc}
a_j & 0\\
0 & 0
\end{array}\right].
$$ 
(Of course, in the beginning of the induction the latter case does not occur.) 
After multiplying $(u,v)$ by a matrix in 
$\Ga\cap Sp(2)\times ...\times Sp(2)$  
we can assume that every matrix 
$$
M_j(\chi)=\left[\begin{array}{cc}
a_j & b_j\\
c_j & d_j
\end{array}\right]= 
\left[\begin{array}{cc}
a_j & 0\\
0 & d_j
\end{array}\right]
$$
is diagonal. We now argue inductively. Suppose that 
$j\in \{2,..,n\}$. We let $d_j':= d_j/gcd(|d_1|,|d_j|)$. Then there are integers 
$\al_j,\be_j$ such that $\al_j d_j' - \be_j d_1'=1$. It follows that 
$$
\left[\begin{array}{cccc}
\al_j & 0 &\be_j & 0\\
a_j & d_j' & -a_1 & -d_1'\end{array}\right] 
\left[\begin{array}{cc}
a_1 & 0\\
0 & d_1\\
a_j & 0\\
0 & d_j
\end{array}\right]= 
\left[\begin{array}{cc}
\al_j a_1 +\be_j a_j & 0\\
0 & 0
\end{array}\right]. 
$$
Note that the row vectors $p,q$ of the first matrix in the above 
formula are such that $\om(p,q)=1$. Hence, according to Lemma 
\ref{completion}, there exists a matrix 
$$
A=\left[\begin{array}{cccc}
* & * & * & *\\
* & * & * & *\\ 
\al & 0 &\be & 0\\
a_2 & d_2' & -a_1 & -d_1'\end{array}\right]
$$
which belongs to $Sp(4,\Z)$. We extend the matrix $A$ to a matrix 
$g\in Sp(2n,\Z)$ which preserves all the coordinates except $a_1,b_1$ 
and $a_j,b_j$. Then the character $\chi'=g\chi$ has $\om_j(\chi')=0$. 
Continuing inductively we find $h\in \Ga$ such that the character
$h\chi$ satisfies:
$$
\om_j(h\chi)=0, j=2,3,...,n. 
$$ 
Note that $\om_1(h\chi)= \om(h\chi)=\om(\chi)>0$. 
Recall that $\Image(\chi)=\Z\times\Z$. Hence $b_1'=\chi(y_1)=1$, since 
 all other generators $x_1, x_2,y_2,...$ of $H_1(S,\Z)$ are mapped by 
$h\chi$ to the real numbers.  
Finally, to obtain $\ga\chi$ as required by the lemma, we multiply $h\chi$ 
by a diagonal symplectic matrix with diagonal entries in $\{\pm 1\}$ to get 
$a_j\ge 0$ for $j=2,...,n$. \qed 

We again use the notation $\chi$ for the character $\chi'$ obtained 
in the previous lemma. 

\begin{lem}
There exists $\ga\in\Ga$ such that that the character 
$\ga\chi$ satisfies:

(i) 
$$
M_1(\ga\chi)= \left[\begin{array}{cc}
a_1=\om(\chi) & 0\\
0 & 1
\end{array}\right].
$$
(ii) For each $j\ge 2$, 
$$
M_j(\ga\chi')=  
\left[\begin{array}{cc}
a_j' & 0\\
0 & 0
\end{array}\right], 0\le a_j' < a_1. 
$$
\end{lem}
\proof For each $j\ge 2$ there exists $t_j\in \Z$ such that 
$0\le a_j':=a_j- t_j a_1 <a_1$.  
Then form the symplectic matrix
$$
\ga=\left[\begin{array}{cc|cc|cc|c|cc}
1 &   0 & 0 &  0 &   0 &  0   &   \ldots &  0 &  0\\
0 &   1 & 0 & t_2  & 0 &  t_3 &    \ldots &  0 &  t_n\\
-t_2 &0 & 1 &  0 &   0 &   0&     \ldots &  0 &  0\\
0 &   0 & 0 &  1 &   0 &   0&     \ldots &  0 &  0\\
-t_3 &0 & 0 &  0 &   1 &   0&    \ldots & 0 &  0\\
0 &   0 & 0 &  0 &   0 &   1 &    \ldots & 0 &  0\\
\vdots & \vdots & \vdots & \vdots & \vdots & \vdots & \vdots & \vdots & \vdots \\ 
-t_n & 0 & 0 & 0 &   0 &  0&     \ldots  & 1 & 0\\
0 &    0 & 0 & 0 &   0 &   0 &     \ldots & 0 & 1\\
\end{array}\right]. 
$$
The reader will note that this matrix belongs to the Heisenberg subgroup 
of $Sp(2n)$ associated to the flag $(\Span(e_1,e_2), \Span(e_2))$. Then 
$\ga\chi$ has the requires properties:
$$
\left[\begin{array}{cccc}
1 & 0 &0 & 0\\
0 & 1 & 0 & t_j\\
-t_j & 0 & 1 & 0\\
0 & 0 & 0 & 1
\end{array}\right] 
\left[\begin{array}{cc}
a_1 & 0\\ 0 & 1\\ 
a_j & 0 \\ 0 & 0 
\end{array}\right]
=
\left[\begin{array}{cc}
a_1 & 0\\ 0 & 1\\ 
a_j' & 0 \\ 0 & 0 
\end{array}\right]. \qed 
$$
We note that for some $j$ we might have $a_j'=0$. However, since  
$\om(\chi)\ge 2=\Area(\C/\Z^2)$ we conclude that there exists 
at least one $j\ge 2$ such that $a_j>0$. Rename this index $j$ to make it 
equal to $2$. Rename $\chi'=\ga\chi$ back to $\chi$  and $a_j'$ back 
to $a_j$, $j=2,...,n$. 

\begin{lem}
There exists $\ga\in\Ga$ such that that the character 
$\ga\chi$ satisfies: 

(i) 
$$
M_1(\ga\chi)= \left[\begin{array}{cc}
a_1=\om(\chi) & 0\\
0 & 1
\end{array}\right].
$$

(ii) For each $j\ge 2$, 
$$
M_j(\ga\chi)=  
\left[\begin{array}{cc}
a_j' & 0\\
0 & 0
\end{array}\right], \quad 0< a_j' <a_1. 
$$
\end{lem}
\proof  The required matrix $\ga$ belongs to the Heisenberg group 
associated to the flag  $(\Span(e_3,e_4), \Span(e_4))$. For each $j$ such 
that $a_j\ne 0$ the multiplication by $\ga$ will not change $a_j$ at all. 
Suppose that $j\ge 3$, $a_j=0$. 
We describe the case $j=3$ and $n=3$, the general case is done inductively.  
$$
\ga= \left[\begin{array}{cccccc}
1 & 0 & 0 & 0 & 0 & 0\\
0 & 1 & 0 & 0 & 0 & 0\\
0 & 0 & 1 & 0  &0 & 0\\
0 & 0 & 0 & 1  &0 & -1\\
0 & 0 & 1 & 0 & 1 & 0\\
0 & 0 & 0 & 0 & 0 & 1
\end{array}\right].$$
Then
$$
\ga M(\chi)= 
\left[\begin{array}{cc}
a_1 & 0 \\ 0 & 1 \\ 
a_2 & 0\\ 0 & 0 \\
a_2 & 0\\ 0 & 0 
\end{array}\right]. \qed 
$$

\section{The non-semisimple case}\label{reduciblecase}

In this section we analyze lattices in those non-semisimple Lie subgroups 
$F$ of $Sp(2n,\R)$ that contain $Sp(2n-2,\R)$. Recall that each maximal 
non-semisimple subgroup $B$ of $Sp(2n,\R)$ containing $Sp(2n-2,\R)$, 
preserves a line $L\subset V^{\perp}$, where $V=\R^{2n-2}$ is the 
symplectic subspace invariant under $Sp(2n-2)$. The group $B$ splits 
as semi-direct product $CH_{2n}\rtimes Sp(2n-2)$, where 
$CH_{2n}=H_{2n-1}\rtimes \R_+$ and $H_{2n-1}$ is the $2n-1$-dimensional 
Heisenberg group, see \S \ref{TH}. 

Now suppose that $F=F_\chi\subset B$ is a Lie subgroup containing $Sp(2n-2)$. 
Since $Sp(2n-2)$ acts transitively on $V-0$, 
the subgroup $F$ has to be one of the following:  

(a) $F=B$. 

(b) $F= H_{2n-1}\rtimes Sp(2n-2)$. 

(c) $F= A\times Sp(2n-2)$ where $A\subset B_0=\R\rtimes \R_+$. 

\medskip
If $\Del=F\cap Sp(2n,\Z) \subset F$ is a lattice then its intersection 
with the subgroup 
$CH_{2n}$ (case (a)), $H_{2n-1}$ (case (b)) and $A$ (case (c)) is again a 
lattice (see Theorem \ref{rw}). The first case is impossible by 
Lemma \ref{flow}. In the third case 
the intersection $\Del\cap Sp(2n-2)$ is a lattice as well and we are therefore 
reduced to the discussion in \S \ref{reduciblecase}. 
This leaves us with the case (b), when $Sp(2n,\Z)\cap H_{2n-1}$ is a lattice. 
Note that there are lattices $\Del\subset  H_{2n-1}\rtimes Sp(2n-2)$ 
whose intersection with any conjugate of $Sp(2n-2)$ is not a lattice, 
we leave it to the reader to construct such examples.

Suppose now that $\chi\in X$ is a character 
(with the real part $u$ and the imaginary part $v$) such that the closure of 
the orbit $\Ga\chi$ contains the orbit $F_\chi \chi$ where 
$F_\chi\cong  H_{2n-1}\rtimes Sp(2n-2)$ fixes a line $L$ in $\Span(u,v)$. 
According to  Remark \ref{rem1} it suffices to consider the case 
$L=\Span(u)$. Applying an element $\ga\in \Ga$ we can adjust the pair $(u,v)$ 
such that the vector $u=(a_1,b_1,...,a_n,b_n)$ 
has no zero coordinates (see Lemma \ref{zero}). The group $H_{2n-1}$ acts 
transitively on the set of vectors $v\in \R^{2n}$ satisfying $\om(u,v)=1$. 
Hence we can find $h\in H_{2n-1}$ such that 
$$
h(v)= \frac{1}{\om(u,v)}(....,-b_j, a_j,....).  
$$
Hence $\om_j(u,h(v))=\om_j(h(u),h(v))>0$ 
for each $j=1,...,n$. Since $\ol{\Ga\chi}$ contains the orbit 
$F_\chi \chi$, there exists an element $\ga\in\Ga$ such that 
$\om_j(\ga(u),\ga(v))>0$ for each $j=1,...,n$.
According to Theorem \ref{geo} 
the character $\chi$ belongs to the subset $\Si\subset X$ of 
characters of abelian differentials.

\section{Meromorphic differentials}\label{mero}

\begin{thm}\label{pole}
Suppose that $n\ge 3$ and $\chi$ is a nonzero character in $H^1(S,\C)$ which does not 
satisfy either Obstruction 1 or Obstruction 2. Then there is a complex 
structure $\tau$ on $S$ and a meromorphic differential $\al$ with a single 
simple pole on $S_\tau$ such that $\chi$ is the character of $\al$.  
\end{thm}
\proof 
Case A. The vectors $u$ and $v$ are linearly independent. 
The group $Sp(2n,\R)$ acts transitively on 
the collection $Y$ of pairs of vectors $u,v\in \R^{2n}$ such that $\om(u,v)=0$ 
and $u\wedge v\ne 0$. 
 Thus (since $\Ga=Sp(2n,\Z)$ is Zariski dense in $Sp(2n,\R)$)  there exists 
$\ga\in\Ga$ such that $\chi'=\ga\chi$ satisfies: 
$\om_j(\chi')\ne 0$ for each $j=1,...,n$. 
If each $\om_j(\chi')>0$ then $\chi$ is the character of an abelian 
differential and there is nothing to prove. Hence (after relabelling $j$'s) 
we get: $\om_1(\chi')<0$ and $\om_j(\chi')\ne 0$, $j=2,...,n$. Set 
$\chi:=\chi'$.

\begin{figure}[tbh]
\centerline{\epsfxsize=5in \epsfbox{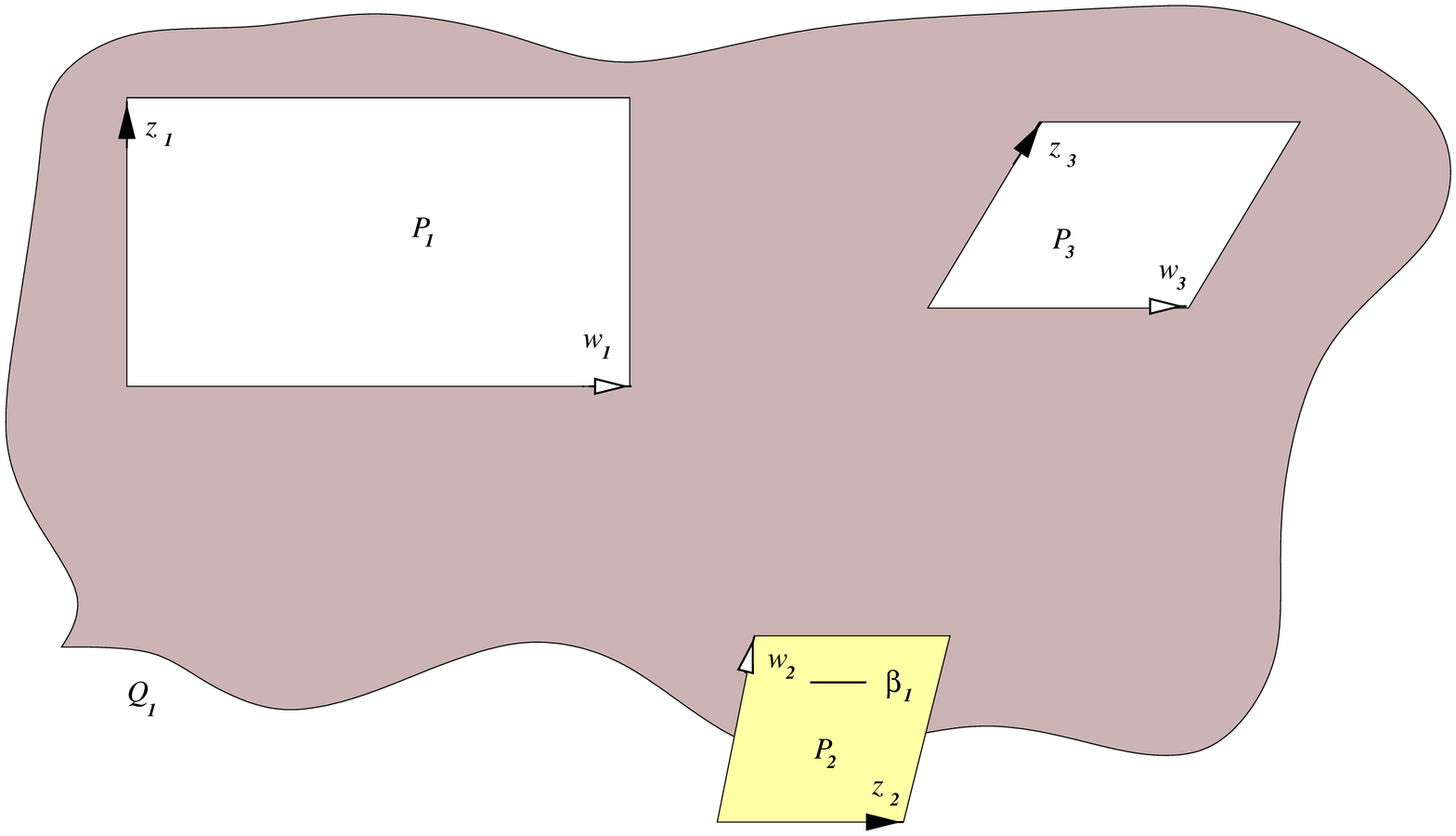}}
\caption{\sl }
\label{F3}
\end{figure}

We argue similarly to the proof of 
Theorem \ref{geo}. Consider the fundamental parallelogram $P_1\subset \C$  
for the discrete group generated by the columns $z_1,w_1$ 
of the matrix $M_1(\chi')$. Let $Q_1$ denote the closure of the 
{\em exterior} of $P_1$ in ${\mathbb S}^2$. Note that topologically 
$Q_1$ is still a parallelogram: its edges are the edges of $P_1$. 
Identifying the opposite sides of $Q_1$ by $z_1,w_1$ we get a marked 
torus $T_1$ with a standard (symplectic) system of generators $x_1,y_1$,  
branched projective structure and an orientation-preserving 
 developing mapping to ${\mathbb S}^2$ whose holonomy is the homomorphism 
$\chi_1$ which sends $x_1\to z_1, y_1\to w_1$. 
(Here we identify a vector in $\C$ with the 
corresponding translation.) Taking pull-back of the form $dz$ on $\C$ we get 
a meromorphic differential on $T_1$ with the single simple pole (corresponding 
to the point $\infty\in Q_1$) and the period character $\chi_1$. We now 
extend this to the rest of the surface $S$. If $j\ge 2$ is such that 
$\om_j(\chi)>0$ then similarly to the proof of Theorem \ref{geo} we add to 
$T_1$ the flat torus $T_j$ obtained by identifying the sides of a fundamental 
parallelogram for the translation group generated by the columns of 
$M_j(\chi)$. If $\om_j(\chi)<0$ we pick a fundamental parallelogram $P_j$ so 
that it is disjoint from the $P_i$'s ($1\le i\le n$, $i\ne j$). Remove 
the interior of $P_j$ from $Q_1$ and identify the opposite 
sides of $P_j$ via translations. See Figure \ref{F3}.

As the outcome we get an oriented surface $S$ and  
a $\chi$-equivariant developing map to ${\mathbb S}^2$. 
The meromorphic differential on $S$ is obtained via pull-back of 
$dz$ from $\C$. Its only 
pole corresponds to the point on the torus $T_1$ which maps to $\infty$ 
under the developing map.

\begin{figure}[tbh]
\centerline{\epsfxsize=5in \epsfbox{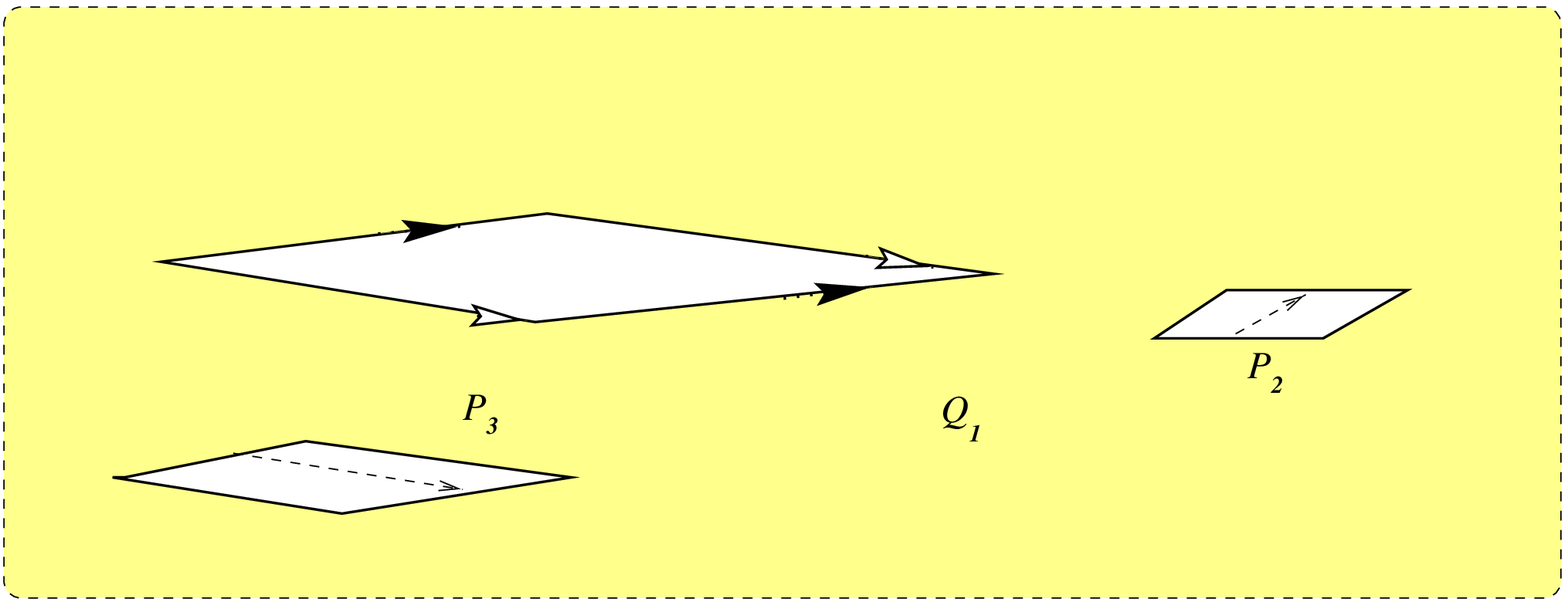}}
\caption{\sl }
\label{F4}
\end{figure}

Case B.  Let $u$ and $v$ be linearly dependent. It suffices to consider 
the case $u\ne 0$ (otherwise replace $\chi$ by $\sqrt{-1}\chi$).  
Using Zariski density of $\Ga$ in $Sp(2n,\R)$ 
(the latter acts transitively on $\R^{2n}-0$) choose $\ga\in\Ga$ 
such that no coordinate of $\ga(u)$ is zero and let $\chi:= \ga\chi$. We now 
argue  analogously to the Case A. Let $z_1,w_2$ denote the columns of the 
matrix $M_1(\chi)$.  Let $P_1$ denote the 
convex hull of the set $0, z_1, w_1, z_1+w_1$. 
  We will think of $P_1$ as a degenerate parallelogram with the edges 
$[0,z_1], [0,w_1], [z_1, z_1+w_1], [w_1, z_1+w_1]$. Now  
cut ${\mathbb S}^2$ open along $P_1$ and denote the result $Q_1$, it is 
homeomorphic to a parallelogram, identification of the opposite edges via 
translations by $z_1, w_1$ yields the torus $T_1$. To reconstruct the rest of 
the surface $S$ we choose disjoint degenerate ``fundamental parallelograms'' 
$P_j$ for the groups generated by the translations $z_j, w_j$, cut $Q_1$ 
open along the $P_j$'s ($j\ge 2$) and get $S$ by identifying the opposite 
edges on each cut.  See Figure \ref{F4}. \qed

\begin{rem}
We note that the branched projective structures $\si$ associated to the 
meromorphic differentials constructed in the above theorem have the branching 
degree $\deg(\si)=2n$. 
\end{rem}

\medskip
We will next prove a lower bound on the degree of branching of the 
projective structures with the holonomy in the translation subgroup $\C$ 
of $PSL(2,\C)$. This lower bound holds for all genera $n\ge 2$. 

Suppose that $\si$  is a branched projective structure 
with the holonomy $\rho: \pi_1(S)\to \C\subset PSL(2,\C)$. We will assume 
that $\rho$ is nontrivial, otherwise clearly $\deg(\si)\ge 2n+2$ by the 
Riemann--Hurwitz formula. The representation $\rho$ lifts to a 
representation $\theta: \pi_1(S)\to SL(2,\C)$ (with the image in 
the group of unipotent upper triangular matrices $U$).  
Let $V$ denote the holomorphic $\C^2$--bundle over $S$ associated with 
the representation $\theta$. The structure $\si$ gives rise to 
a holomorphic line subbundle $L\subset V$ such that
\begin{equation}\label{bun}
\deg(L)= n-1 - \frac{\deg(\si)}{2}, 
\end{equation}
where $\deg(\si)$ is the degree of branching of $\si$ 
(see \cite[Chapter C]{GKM}). The bundle $V$ fits into short exact sequence 
$$
0 \to \La \to V \stackrel{p}{\to} \La \to 0, 
$$
where $\La$ is the trivial bundle; the fibers of 
$\La=\ker(p)$ correspond to the line in $\C$ fixed 
by the group $U$. Under the projectivization $\C^2\to \C\P^1$ 
this line projects to the point $\infty\in  \C\P^1$. Hence the developing 
mapping of $\si$ does not cover $\infty$ iff $L\cap \ker(p)=0$. It also 
follows that $L\ne \ker(p)$ (otherwise the developing mapping of $\si$ would 
be constant). Therefore we get a nonzero map $p: L\to \La$ by restricting 
the projection $p: V\to \La$ to $L$. By the Riemann--Roch theorem,  $\deg(L)\le 0$ with the 
equality iff $p: L\to \La$ is injective; (\ref{bun}) then implies that 
$\deg(\si)\ge 2n-2$. The equality here is attained only if the developing map 
of $\si$ takes values in $\C$, i.e. $\si$ is a singular Euclidean structure. 
In other words, if $\deg(\si)=2n-2$ then the developing mapping 
of $\si$ is obtained by integrating an abelian differential on $S$. 
If $\rho$ is not the holonomy of any singular Euclidean structure 
then $\deg(\si)\ge 2n+1$. However, since $\rho$ lifts to $SL(2,\C)$, 
$\deg(\si)$  has to be even (see  \cite[Chapter C]{GKM}). We conclude that 
in this case $\deg(\si)\ge 2n$. Recall that for a representation 
$\rho: \pi_1(S) \to PSL(2,\C)$, $d(\rho)$ is the least degree of branching 
of all projective structures on $S$ (consistent with the orientation) 
with the holonomy $\rho$. We thus proved:

\begin{prop}
Suppose that $\rho$ is a representation $\rho: \pi_1(S) \to PSL(2,\C)$  
whose image is contained the translation subgroup $\C$ 
of $PSL(2,\C)$. Then 
$d(\rho)\ge 2n-2$ and $d(\rho)\ge 2n$ provided that the corresponding 
character $\chi\in H^1(S,\C)$ is not the period character of any abelian 
differential. 
\end{prop}

Combining this proposition with Theorems \ref{main} and \ref{pole} we obtain   
Corollary \ref{corfor}.

\bibliography{lit}

\bibliographystyle{amsalpha}

\vspace*{2in}

\vspace{3\baselineskip}\noindent
Department of Mathematics,
University of California,
1 Shields Ave., Davis, CA 95616,
USA\\

\vspace*{0.3in}

\noindent kapovich@math.ucdavis.edu

\end{document}